\input amstex

\def\lr#1#2{{\vcenter{\vbox{\hbox{\vrule width0pt height#1pt \kern#1pt
  \vrule width.#2pt}\hrule height.#2pt}}}}
\def\intprod{\mathbin{\lr54}}
\define\worm{\Cal W}
\redefine\Re{\operatorname{Re}}
\redefine\Im{\operatorname{Im}}
\define\sgn{\operatorname{sgn}}

\define\k#1.{\ref\key{#1}}
\define\jj#1.{\langle{#1}\rangle}
\define\half{ {1/2} }
\define\er{\endref}

\define\CC{\bold C}
\define\p{\partial}
\define\nf{\infty}
\define\rp{ {^{-1}} }
\define\al{\alpha}
\define\be{\beta}
\define\de{\delta}
\define\e{\varepsilon}
\define\g{\gamma}
\define\z{\zeta}
\define\la{\lambda}
\define\Om{\Omega}
\define\om{\omega}

\define\lan{\langle}
\define\th{\theta}
\define\ra{\rangle}
\define\lt{ {L^2} }

\define\bl{\Bbb L}
\define\dt{ t^{-1}dt }
\define\okay{\frak B}
\define\jz{\langle\z\rangle}
\define\rr{\bold R}

\documentstyle{amsppt}
\refstyle{A}
\widestnumber\key{ABC}
\magnification=\magstep1
\pageheight{8.6truein}
\pagewidth{6.5truein}
\NoRunningHeads
\NoBlackBoxes

\topmatter
\title Global $C^\nf$ Irregularity of the 
$\bar\partial$--Neumann Problem For Worm Domains
\endtitle
\author Michael Christ\endauthor
\affil University of California, Los Angeles\endaffil
\address Department of Mathematics, UCLA, 
Los Angeles, Ca. 90095-1555 \endaddress
\email christ \@ math.ucla.edu\endemail
\thanks Research supported by National
Science Foundation grant DMS-9306833. I am indebted to D.~Barrett,
J.~J.~Kohn, P.~Matheos and J.~McNeal for stimulating conversations
and useful comments on the exposition. 
\endthanks
\endtopmatter

\document
\subheading{0. Introduction}

Let $\Om\subset\CC^n$ be a bounded, pseudoconvex domain with
$C^\nf$ boundary, equipped with the standard 
Hermitian metric inherited from $\CC^n$.
The $\bar\p$--Neumann problem for $(p,q)$ forms in $\Om$
is the boundary value problem
$$ \left\{ \alignedat 2 
\square u &=f \qquad &&\text{in }\Om  \\
u \intprod \bar\p\rho&=0 \qquad &&\text{on }\p\Om  \\
\bar\p u\intprod \bar\p\rho&=0 \qquad &&\text{on }\p\Om 
\endalignedat 
\right.
$$
where $\rho$ is a defining function for $\Om$,
$\square = \bar\p\bar\p^*+ \bar\p^*\bar\p$, $u,f$ are $(p,q)$ forms, 
and $\intprod$ denotes the interior product of forms.
Under the stated hypotheses on $\Om$, this problem is uniquely
solvable for every $f\in\lt(\Om)$. The Neumann operator $N$,
mapping $f$ to the solution $u$, is continuous on $\lt(\Om)$.
The Bergman projection $B$ is the orthogonal projection of
$\lt(\Om)$ onto the closed subspace of $\lt$ holomorphic functions
on $\Om$, and is related to $N$ by $B = I - \bar\p^* N\bar\p$.

$N$ and $B$ are $C^\nf$ pseudolocal if
$\Om$ is strictly pseudoconvex, or more generally, is of finite type 
\cite{Ca1}. Both preserve 
$C^\nf(\overline{\Om})$ under certain weaker hypotheses \cite{BS2},\cite{Ca2}.
For any pseudoconvex, smoothly bounded $\Om$ 
and any finite exponent $s$, there exists a strictly positive weight 
$w\in C^\nf(\overline{\Om})$ 
such that the Neumann operator and Bergman projection with respect
to the Hilbert space $\lt(\Om,w(x)dx)$ map the Sobolev space
$H^t(\Om)$ boundedly
to $H^t(\Om)$, for all $0\le t\le s$ \cite{K1}.  
It has remained an open question whether $N$ and $B$, 
defined with respect to the standard metric, preserve 
$C^\nf(\overline{\Om})$ without further hypotheses on $\Om$.
An affirmative answer would have significant consequences \cite{BL}.

\proclaim{Theorem} There exist pseudoconvex,
smoothly bounded domains $\Om\Subset\CC^2$ 
for which the Neumann operator on $(0,1)$ forms 
and Bergman projection fail to preserve $C^\nf(\overline{\Om})$.
\endproclaim

Examples are the worm domains, originally 
introduced by Diederich and Forn\ae ss \cite{DF}
for another purpose\footnote{Some but not all 
worm domains have nontrivial Nebenh\"ulle\cite{FS, p\. 111},
whereas all worm
domains are counterexamples to global regularity.} 
but long considered likely 
candidates for $N$ and $B$ to fail to be globally regular in $C^\nf$. 

The proof depends on the observation of
Barrett \cite{B} that for each worm domain $\Cal W$,
for all sufficiently large $s$, $N$ and $B$ fail
to map $H^s(\Cal W)$ boundedly to $H^s(\Cal W)$.
We establish for each worm domain
an {\it a priori} inequality
of the form $\|Nf\|_{H^s} \le C_s\|f\|_{H^s}$, valid for
all $f\in C^\nf(\overline{\Cal W})$ such that 
$Nf\in C^\nf(\overline{\Cal W})$,
for a sequence of exponents $s$ tending to $\nf$. If $N$ were
to preserve $C^\nf(\overline{\Cal W})$, then since it is a bounded
linear operator on $\lt(\Cal W)$ and since $C^\nf(\overline{\Cal W})$
is dense in $H^s(\Cal W)$, it would follow that $N$
maps $H^s(\Cal W)$ boundedly to itself, for a sequence
of values of $s$ tending to $\nf$, a contradiction.
More accurately, our inequality is valid only for certain subspaces of 
$\lt(\worm)$ preserved by $\square$, but this still suffices to 
contradict the theorem of Barrett.

An analogous counterexample in the real analytic context
is already known \cite{Ch1}: there exists a bounded, pseudoconvex
domain $\Om\subset\CC^2$ having real analytic boundary,
such that the Szeg\"o projection fails to preserve
$C^\om(\p\Om)$. That result and its proof are however
not closely related to the $C^\infty$ case.

\S\S 1 through 3 review material on worm domains 
and the $\bar\partial$--Neumann problem, and present some
routine but tedious reductions. 
\S4 formalizes a class of two-dimensional problems
subsuming those to which the reductions lead.
The analysis of those problems is contained in \S\S 5 and 6. 

\subheading{1. Reduction to the Boundary}

The $\bar\p$--Neumann problem is a boundary value problem
for an elliptic partial differential equation, and as such
is amenable to treatment by the method of reduction to
a pseudodifferential equation on the boundary.
This reduction has been carried out in detail for domains
in $\CC^2$ by Chang, Nagel and Stein \cite{CNS}. We review here
certain of their computations and direct consequences thereof.

Assume $\Om\subset\CC^2$ to be a smoothly bounded domain.
The equation $\square u=f$ on $\Omega$ for $(0,1)$ forms in
$C^\nf(\overline{\Om})$ is equivalent to an equation 
$\square^+ v=g$ on $\p\Om$, where $v,g$ are sections of
a certain complex line bundle\footnote{$\Cal B^{0,1}$
is defined to be the quotient of the
restriction to $\p\Om$ of $T^{0,1}\CC^2$,
modulo the span of $\bar\p\rho$. Sections of $\Cal B^{0,1}$
may be identified with scalar-valued functions times $\bar\omega_1$,
hence with scalar-valued functions.}
$\Cal B^{0,1}$. Let $\rho$ be a 
smooth defining function for $\Om$ and define 
$\bar\omega_2 = \bar\p\rho$, 
and $\bar\omega_1 = (\p\rho/\p  z_2)d\bar z_1
- (\p\rho/\p  z_1)d\bar z_2$.  

$v$ is related to $u$ by $u=Pv+Gf$, where $P,G$ are
respectively Poisson and Green operators for the elliptic
system $\square u=f$ with Dirichlet boundary conditions. 
In particular, if $f\in C^\nf(\overline{\Om})$, then
$u\in C^\nf(\overline{\Om})$ if and only if the same holds for $v$.
More precisely,
$G$ maps $H^s(\Om)$ to $H^{s+2}(\Om)$,
while $P$ maps $H^s(\p\Om)$ to $H^{s+\half}(\Om)$,
for each $s\ge 0$.
Thus if $f\in H^s$, in order to conclude that $u\in H^s$
it suffices to know that $v\in H^{s-\half}$.

On the other hand, $g = (\bar\p Gf\intprod\bar\p\rho)$,
restricted to $\p\Om$. If $f\in H^s$ then $\bar\p Gf\in H^{s+1}(\Om)$,
so its restriction to the boundary belongs to $H^{s+\half}(\p\Om)$.
Thus in order to show that $N$ preserves $H^s(\Om)$ it suffices to
show that if $\square^+ v\in H^{s+\half}(\p\Om)$, then $v\in
H^{s-\half}(\p\Om)$, assuming always that $s>\half$.

On $\p\Om$ a Cauchy-Riemann operator is the complex vector field 
$\bar L = 
(\p_{\bar z_1}\rho) \p_{\bar z_2} 
- (\p_{\bar z_2}\rho) \p_{\bar z_1} $.
Define $L$ to be the complex conjugate of $\bar L$.
The characteristic variety \footnote{By the characteristic
variety of a pseudodifferential operator we mean the conic
subset of the cotangent bundle on which its principal symbol
vanishes.}
of $\bar L$ is a real line bundle $\Gamma$. 
Assuming $\Om$ to be pseudoconvex and the set of strictly pseudoconvex
points to be dense in $\p\Om$,
$\Gamma$ splits smoothly and uniquely as $\Gamma^+\cup\Gamma^-$,
where each fiber of $\Gamma^\pm$ is a single ray, 
and where $\Gamma^+$ is distinguished from $\Gamma^-$ by the requirement
that the principal symbol of $[\bar L, L]$ is nonpositive on $\Gamma^+$,
modulo terms spanned by the symbols of the real and imaginary
parts of $\bar L$ and a term of order $0$.
Equivalently, the principal symbol of $[\bar L, \bar L^*]$
is nonnegative on $\Gamma^+$, modulo the same kinds of
error terms.

$\square^+$ is a classical pseudodifferential operator of
order $+1$. Its principal symbol vanishes everywhere on $\Gamma^+$
but nowhere else. Microlocally in a conic neighborhood of
$\Gamma^+$, $\square^+$ takes the form
$$
\square^+ = Q\bar L L + F_1 \bar L + F_2 L + F_3
$$
where $Q$ is an elliptic pseudodifferential operator of order $-1$,
and each $F_j$ is a pseudodifferential operator of order less than
or equal to $-1$.
Since $\square^+$ is elliptic except on $\Gamma^+$,
for any pseudodifferential operator $G$ of order zero whose symbol
vanishes identically in some neighborhood of $\Gamma^+$,
one has for all $u\in C^\nf$ and all $N<\nf$
$$
\|Gu\|_{H^{t+1}(\p\Om)} \le C \|\square^+ u\|_{H^t(\p\Om)}
+ C_N\|u\|_{H^{-N}(\p\Om)}.
\tag1.1
$$

Let $A$ be an elliptic pseudodifferential operator of order $+1$
such that $A\circ Q$ equals the identity on $\lt(\p\Om)$, modulo
an operator smoothing of infinite order. Composing on the left with
$A$, the equation $\square^+ v=g$ may be rewritten as
$\frak L v= \tilde g$ microlocally in a conic neighborhood of
$\Gamma^+$, where 
$$
\frak L = \bar L L + B_1\bar L + B_2 L + B_3,
\tag1.2
$$
$\|\tilde g\|_{H^t} \le C \|g\|_{H^{t+1}} + C_N\|v\|_{H^{-N}}$ 
for any finite $N$,
and each $B_j$ is an operator of order less than or equal to zero.
Therefore in order to show that the Neumann operator satisfies
an {\it a priori} inequality 
of the form $\|Nf\|_{H^s(\Om)} \le C\|f\|_{H^s(\Om)}$ for
all $f\in C^\nf(\overline{\Om})$ such that 
$Nf\in C^\nf(\overline{\Om})$,
it suffices to establish an {\it priori inequality}
for all $v\in C^\nf(\p\Om)$ of the form
$$
\|v\|_{H^t}
\le C\|\frak L v\|_{H^t} + C \|v\|_{H^{t'}}
+ C\|\tilde Qv \|_{H^{t+2}}
\tag1.3
$$
where $t=s-\half$, 
for some $t'<t$ and some pseudodifferential operator
$\tilde Q$ of order zero whose symbol vanishes identically
in some neighborhood of $\Gamma^+$.

\subheading{2. Worm Domains}

A worm domain in $\CC^2$ is an open set of the form
$$
\Cal W =\{z: |z_1 + e^{i\log|z_2|^2}|^2 < 1-\phi(\log|z_2|^2)\}
$$
where the function $\phi$ vanishes identically on some
interval $[-r,r]$ of positive length, and is constructed \cite{DF}
so as to guarantee that $\Cal W$ will be pseudoconvex with $C^\nf$
boundary, and will be strictly pseudoconvex at every boundary
point except those on the exceptional annulus $\Cal A\subset\p\Cal W$
defined as
$$
\Cal A = \{z: z_1=0\text{ and } |\log |z_2|^2 |\le r\}.
$$

The circle group acts as a group of automorphisms of $\Cal W$
by $z\mapsto R_\th z = (z_1,e^{i\theta} z_2)$.
It acts on functions by $R_\th f(z) = f(R_\th z)$, and
on $(0,1)$ forms by $R_\th(f_1 d\bar z_1 + f_2 d\bar z_2)
= (R_\th f_1)d\bar z_1 + (R_\th f_2)e^{-i\th} d\bar z_2$.
The Hilbert space $L^2_{(0,k)}(\worm)$ of square integrable 
$(0,k)$ forms
decomposes as the orthogonal
direct sum $\oplus_{j\in\bold Z}\Cal H_j^k$ 
where $\Cal H_j^k$ is the set of all $(0,k)$ forms $f$ 
satisfying $R_\th f \equiv e^{ij\th}f$.
$\bar\p$ is an unbounded linear operator from $\Cal H_j^k$ to
$\Cal H_j^{k+1}$, $B$ maps $\Cal H_j^0$ to itself,
and the Neumann operator $N$ maps
$\Cal H_j^1$ to $\Cal H_j^1$ boundedly, for each $j$. 

For each $k$ and each $s\ge 0$, the Sobolev space
$H^s(\worm)$ likewise decomposes as an orthogonal direct
sum of subspaces $H^s_j$.
It is known that for any smoothly bounded, pseudoconvex domain
$\Om\subset\CC^2$, for any exponent $s\ge 0$, if $N$ maps
$H^s(\Om)$ boundedly to itself, then $B$ also maps
$H^s(\Om)$ boundedly to itself \cite{BS1}, where $H^s$ denotes 
in the first instance a space of one forms, and in the 
second, a space of functions.
Because $N,B$ preserve the summands $\Cal H_j$, the same proof
shows\footnote{This follows from the argument of Boas and Straube
\cite{BS1} because all elements of their proof may
be chosen to be invariant under the automorphisms $R_\th$.}
that for any fixed $j$, if $N$ maps the space $H^s_j$ of
$(0,1)$ forms boundedly to itself, then $B$ maps the space
$H^s_j$ of functions boundedly to itself. Barrett \cite{B}
that for each worm domain, for all sufficiently large
$s$, for all $j$,
$B$ fails to map $H^s_j$ boundedly to itself. Therefore
in order to prove that $N$, acting on $(0,1)$ forms, fails to 
preserve $C^\nf(\overline{\worm})$, 
it suffices to establish the following result for a single index $j$.

\proclaim{Proposition 1}
For each worm domain 
there exists a discrete subset $S\subset\rr^+$
such that for each $s\notin S$ 
and each $j\in\bold Z$
there exists $C_{s,j}<\nf$
such that for every $(0,1)$ form 
$u\in \Cal H_j^1\cap C^\nf(\overline{\worm})$
such that $Nu\in C^\nf(\overline{\worm})$,
$$
\|Nu\|_{H^s(\worm)} \le C_{s,j} \|u\|_{H^s(\worm)}.  \tag2.1
$$
\endproclaim

The defining function $\rho = 1-\phi(\log|z_2|^2)
- |z_1 + e^{i\log|z_2|^2}|^2$  for $\worm$
is invariant under $R_\th$, as is the $(0,1)$
form $\bar\om_2$ defined above. $\bar\omega_1$ satisfies
$R_\alpha \bar\omega_1 = \exp(-i\alpha)\bar\omega_1$ for all
$\alpha$, but
it may also be made invariant by multiplying
it by the function $(z_1, re^{i\th})\mapsto e^{i\th}$,
which is smooth in a neighborhood of $\overline{\worm}$.
We work henceforth with this modified $\bar\omega_1$.

$\square^+$ commutes with $R_\th$ for all $\th$.
Indeed, $\square^+ v = \bar\p Pv \intprod\bar\p\rho$
\cite{CNS}. 
$\square$ commutes with $R_\th$, hence so must $P$.
$\bar\p$ commutes with $R_\th$, and the Hermitian metric
on $\CC^2$ and $\bar\p\rho$ are likewise $R_\th$--invariant.
Thus all ingredients in the above expression for $\square^+$
are invariant, hence so is $\square^+$ itself.

Identify square integrable sections of $\Cal B^{0,1}$ with 
scalar-valued $\lt$ functions as above, and decompose 
$\lt(\p\worm)= \oplus \Cal H_j(\p\worm)$ where 
$\Cal H_j$ is the subspace
of those functions satisfying $R_\th f \equiv e^{ij\th} f$.
Then $\square^+$ maps $\Cal H_j(\p\worm)\cap C^\nf$ to 
$\Cal H_j(\p\worm)$.
We have seen in \S1 that Proposition 1 would be a consequence of the
validity of (1.3) for all $v\in C^\nf(\p\worm)\cap \Cal H_j$,
for $t=s-\half$.

Fix $j$ and assume henceforth that $u$ belongs to $\Cal H_j^1(\worm)$ 
and to $C^\nf$.
Then the associated boundary function $v$ belongs to 
$\Cal H_j \cap C^\nf(\p\worm)$.
Henceforth we work exclusively on the boundary, and
simplify notation by writing simply $\Cal H_j$ rather
than $\Cal H_j(\p\worm)$.

Note that $\bar L,L$ take $\Cal H_j\cap C^\nf$
to $\Cal H_{j+1}$ and to $\Cal H_{j-1}$, respectively.
The operator $A$ introduced after (1.1)
may be constructed to be $R_\th$-invariant, 
for both $\square^+$ and $\bar L\circ L$ are invariant while
$\bar L,L$ are automorphic of certain degrees, so that
averaging the equation
$\square^+ = Q\bar L L + F_1\bar L + F_2 L + F_3$
with respect to $R_\th d\th$ produces an invariant
$Q$ and $F_3$, and operators $F_1,F_2$ 
automorphic of the appropriate degrees.  Thus
$(\bar LL + B_1\bar L + B_2 L + B_3)v\in H^t$
microlocally near $\Gamma^+$,
where $B_1,B_2,B_3$ map
$\Cal H_j$ to $\Cal H_i$ for $i= j-1,j+1,j$ respectively.

Since $\worm$ is strictly pseudoconvex at all points not in $\Cal A$,
it follows as in \cite{K2} that 
on the complement of any neighborhood of $\Cal A$,
the $H^{s+1}$ norm of $v$ is majorized by 
$C\|\square^+ v\|_{H^{t+1}(\p\worm)} + C\|v\|_{H^{-N}}$,
hence by $C\|\frak L v\|_{H^s} + C\|v\|_{H^{-N}}
+ C\|\tilde Qv\|_{H^{s+2}}$ with $\tilde Q$ as in (1.3). 
This estimate is one derivative stronger than that which we seek.
In particular, it now suffices to control the $H^{s}$ norm
of $v$ in an arbitrarily small neighborhood $U$ of $\Cal A$,
and to do so microlocally near $\Gamma^+$.

Fix a $C^\nf$ cutoff function $\varphi$ supported in a small
neighborhood of $\Cal A$ but identically equal to $1$
in a smaller neighborhood, and fix an open set $V$ disjoint
from a neighborhood of $\Cal A$ such that $\nabla \varphi$ is
supported in $V$.
By Leibniz's rule and the pseudolocality of pseudodifferential operators,
the $H^s$ norm of $\frak L(\varphi v - v)$ is majorized by
$C\|v\|_{H^{s+1}(V)} + C\|v\|_{H^{-N}}$ for any $N<\nf$.  Thus 
by replacing $v$ with $\varphi v$ we may reduce matters to the case
where $v$ is supported in an arbitrarily small neighborhood $W$ of 
$\Cal A$. Therefore it suffices to prove the existence of some $W$
and an exponent $s'<s$ such that (1.3) holds (with $t$
replaced by $s$) for all $v\in C_0^\nf$ supported in $W$.
\medskip

In a neighborhood of $\Cal A$ in $\p\worm$ introduce
coordinates $(x,\th,t)$ where 
$$z_2 = e^{x+i\th}, \qquad z_1=e^{i2x}( e^{it}(1-\phi(2x)) -1)$$ 
with $2|x| < r+\de$ and $|t|<\de$ for some small $\de>0$.
In these coordinates $\Cal A = \{t=0,\ |x|\le r/2 \}$.
Setting
$$\gamma(x,t) = 2\big[e^{-it}-1+\phi(2x)-i\phi'(2x)\big]/
\big[1-\phi(2x)\big],$$
the vector field $\bar L = \p_x+i\p_\th + \gamma\p_t$
annihilates both $z_1$ and $z_2$. 
Hence it differs from what was previously denoted as
$\bar L$ by multiplication on the left by a nonvanishing
factor, which may be verified to have the form $b(x,t)\exp(i\th)$.
For $|x|\le r/2$, $\phi(x)\equiv 0$
and consequently $\gamma(x,t)\equiv 2(e^{-it}-1)$. 
Therefore\footnote{
Boas and Straube \cite{BS2} have shown
the $\bar\p$--Neumann problem to be globally
$C^\nf$ hypoelliptic whenever there exists a real vector
field on the boundary that is transverse to the complex tangent
space and has a certain favorable commutation property.
If $\Re\al(0)$ were to vanish then
$\p_t$ would be such a vector field. Thus nonvanishing of
$\Re\al(0)$ is for our
purpose an essential feature of worm domains.}
$$
\bar L = \p_x+i\p_\th + it\al(t)\p_t\qquad\text{where }|x|\le r/2,
\text{ with } \Re\al(0)\ne 0.
\tag2.2
$$
The representation 
$\frak L = \bar L L + B_1\bar L + B_2 L + B_3$
in terms of the new $\bar L$ remains valid
with modified coefficients $B_i$ that now preserve
each $\Cal H_j$, once the operator formerly denoted by 
$\frak L$ is multiplied by $|b|^{-2}$. 

There exists a unique $C^\nf$ real-valued function $\mu$,
independent of $\th$, such that
$[\bar L,L] = i\mu(x,t)\p_t $ modulo the span of the real
and imaginary parts of $\bar L$. More precisely, since
the coefficients of $\p_\th$ in $[\bar L,L]$ and
in $\Re\bar L$ vanish while the coefficient in $\Im\bar L$
is nowhere zero,
$$[\bar L,L] = i\mu(x,t)\p_t + i\nu(x,t) \Re\bar L  \tag2.3$$
for unique real-valued, $C^\nf$ coefficients $\mu,\nu$.
The pseudoconvexity of $\p\worm$ means that
$\mu$ does not change sign. Replacing $t$ by $-t$
if necessary, we may assume that $\mu\ge 0$.

Fix an integer $k$.
We identify functions of $(x,t)\in\rr^2$
with elements of $\Cal H_k$ via the correspondence
$u(x,t)\mapsto u(x,t)e^{ik\theta}$;
$\p_\th$ then becomes multiplication by $ik$.
For the remainder of the paper we work in $\rr^2$.
Define
$$\Cal D=\{(x,t): |x|\le r/2\text{ and }t=0\}.$$
By incorporating $ik$ into the $B_j$
we may further rewrite $\frak L$, when restricted
to $\Cal H_k$, as
$$
\Cal L_0 = \bar\ell\ell + B_1\bar\ell + B_2\ell + B_3  
$$
where $\bar\ell$ is a complex vector field in $\rr^2$ which, 
for $|x|\le r/2$, takes the form
$$
\bar\ell = \p_x + it\al(t)\p_t.
$$ 
Here $\ell$ denotes the conjugate of $\bar\ell$, each $B_j$ 
is a pseudodifferential operator of order $\le 0$ in $\rr^2$,
and $\Re\al(0)\ne 0$.
Note that the commutator of the real and imaginary parts of
$\bar\ell$ is not forced to vanish identically, because $\al$ is not
real-valued.  We have
$[\bar\ell,\ell] = i\mu(x,t)\p_t + i\nu(x,t)\Re\bar\ell$
with the same coefficients as in (2.3).

Letting $(\xi,\tau)$ be Fourier variables dual to $(x,t)$,
define
$\tilde\Gamma=\{(x,t,\xi,\tau): (x,t)\in\Cal D
\text{ and } \xi=0\}$.
$\Cal L_0$ is not elliptic at points
of $\Cal D$, but is elliptic at most points in its complement;
$\tilde\Gamma$ is the intersection of the characteristic
variety of $\Cal L_0$ with $\{(x,t,\xi,\tau): (x,t)\in\Cal D\}$.
Decompose $\tilde\Gamma = \tilde\Gamma^+\cup\tilde\Gamma^-$
where $\tilde\Gamma^+ = \tilde\Gamma\cap\{\tau>0\}$.
Thus the principal symbol of $i\mu\p_t$, namely $-\mu\tau$,
is nonpositive in a conic neighborhood of $\tilde\Gamma^+$.

\subheading{3. Two Pseudodifferential Manipulations}

By an operator we will always mean a classical
pseudodifferential operator, that is, one whose symbol
admits a full asymptotic expansion in homogeneous terms
of integral degrees. $\sigma_j(T)$ denotes the $j$-th order
symbol of $T$ (in the Kohn-Nirenberg calculus), 
always with respect to the fixed coordinate system $(x,t,\xi,\tau)$.
Henceforth we work under the convention that $A,B,E$
denote operators whose orders are
less than or equal to $0,0,-1$ respectively, whose meanings are
permitted to change freely from one occurrence to the next,
even within the same line.
$A$ denotes always an operator having the additional property 
that $\sigma_0(A)(x,t,\xi,\tau)\equiv 0$ 
for all $(x,t)\in\Cal D$.
Any operator of type $E$ may be regarded as one of type $A$.
Two operators are said to agree microlocally in some conic
open set if the full symbol of their difference vanishes
identically there.

Any operator $B$ may be written in the form 
$B=\be(x) + E\circ\ell + A$ microlocally in a conic neighborhood
of $\tilde\Gamma^+$,
where $\be$ denotes both a $C^\infty$ function and the operator
defined by multiplication by that function. Indeed,
$\sigma_0(B)(x,0,0,\tau)$ depends only on $(x,\sgn(\tau))$
and we define $\be(x)$ to be this quantity for $\tau>0$.
Then where $\tau>0$ and $|x|\le r/2$, 
$\sigma_0(B)(x,0,\xi,\tau)$ is divisible
by $\xi = -i\sigma_1(\ell)(x,0,\xi,\tau)$. 
$\sigma_{-1}(E)(x,0,\xi,\tau)$ is then uniquely determined
for such $x,\tau$ by the equation
$\sigma_0(B) = \be(x) + \sigma_1(\ell)\cdot\sigma_{-1}(E)$.
Define $E$ to be any operator of order $-1$ whose principal 
symbol satisfies this equation when
restricted to $(x,t)\in\Cal D$ and to a conic neighborhood
of $\tilde\Gamma^+$.  Then simply define $A = B-\be-E\circ\ell$.

Writing $B_1 = \be_1(x) + E_1\circ\ell + A_1$
and similarly $B_2 = \be_2(x) + E_2\circ\bar\ell + A_2$,
and expressing $E\ell\bar\ell = E\bar\ell\ell$ plus an operator
of order $\le 0$, we obtain
$\Cal L_0 = (I+E)\bar\ell\ell + (\be_1+A)\bar\ell 
+ (\be_2+A)\ell + B$. Composing both sides with a
parametrix for $I+E$ and modifying the definition of
$\Cal L_0$ to include this factor, we have
$\Cal L_0 = \bar\ell\ell + (\be_1+A)\bar\ell
+ (\be_2+A)\ell + B$. Writing finally $B = \be_3(x)
+ E\circ\ell + A$ results in
$$
\Cal L_0 = \bar\ell\ell + (\be_1+A)\bar\ell
+ (\be_2+A)\ell + (\be_3+A)
\tag3.1
$$
where the $\be_j$ are $C^\nf$ functions depending only on $x$.
\hfill\hfill\qed \medskip

We next reduce the question of {\it a priori}
$H^s$ inequalities to $L^2$, simultaneously for all $s$.  
Fix an operator $Q$ of order $0$ that is elliptic
in some conic neighborhood of $\tilde\Gamma^-$, 
whose symbol vanishes identically in some conic neighborhood 
of $\tilde\Gamma^+$.
Fix an exponent $s>0$, for which we seek an {\it a priori} inequality
for all $v\in C^\infty$ of the form
$$\|v\|_{H^s}\le C\|\Cal L_0 v\|_{H^s} + C\|v\|_{H^{s'}}
+ C\|Qv\|_{H^{s+2}}  \tag3.2 $$
for some $s'<s$.
Having such an inequality for a sequence of exponents $s$ tending
to $+\nf$ would imply Proposition 1 by the preceding discussion.
In particular, the $H^{s+2}$
norm of $Qv$ is already under control by virtue of (1.1),
while the $H^0$ norm of $v$ is harmless because the Neumann operator
is bounded on $\lt$, and $\|v\|_{H^{s'}} \le \e \|v\|_{H^s}
+ C_{\e,N}\|v\|_{H^{-N}}$ for any $\e>0$ and $N<\nf$.

Fix a $C^\nf$, strictly positive function $m=m(\xi,\tau)$,
homogeneous of degree $1$ for large $|(\xi,\tau)|$ and identically
equal to $(1+\tau^2)^\half$ in a conic neighborhood of
$\{\xi=0\}$. Define 
$\Lambda^s$ to be the Fourier multiplier operator
on $\rr^2$ with symbol $m(\xi,\tau)^s$. 

Substituting $v=\Lambda^{-s}u$ and $g=\Lambda^{-s}f$, 
estimation of the $H^s$ norm of $v$, modulo a lower order norm,
in terms of that of $g$ is equivalent to estimation
of the $H^0$ norm of $u$ in terms of that of $f$,
modulo a negative order norm of $u$.
The equation $\Cal L_0 v=g$ becomes
$\Lambda^{-s}\Cal L_0 \Lambda^s u=f$. 
Write $\Lambda^{-s}\bar\ell\ell\Lambda^s =
\Lambda^{-s}\bar\ell\Lambda^s\circ\Lambda^{-s}\ell\Lambda^s$
and similarly for other terms, and note that
$ \Lambda^{-s}(\be_j+A)\Lambda^s
= \be_j+A$ with the same function $\be_j$.

$\Lambda^{\pm s}$ commutes with $\p_x$ and with $\p_t$, so
$ \Lambda^{-s}(\p_x + i\al t\p_t)\Lambda^s =
\p_x + \Lambda^{-s}t\Lambda^s\circ\Lambda^{-s}\al\Lambda^s\p_t $.
Applying the Fourier transform gives immediately
$ \Lambda^{-s}[t,\Lambda^s]\p_t = -s + E $,
microlocally in some conic neighborhood of $\tilde\Gamma^+$,
modulo operators smoothing there of infinite order.
Since
$\Lambda^{-s}[\al,\Lambda^s]\p_t$ is of order $0$, 
$$\align
\Lambda^{-s}t\Lambda^s\circ\Lambda^{-s}\al\Lambda^s\p_t
& = 
\al t\p_t + \Lambda^{-s}[t,\Lambda^s]\al\p_t 
+ t\circ \Lambda^{-s}[\al,\Lambda^s]\circ\p_t
+  E  \\
& =
\al\circ(t\p_t -s) + tB + E,
\endalign
$$ 
microlocally near $\tilde\Gamma^+$.
Thus microlocally near $\tilde\Gamma^+$, 
$\Lambda^{-s}\Cal L_0\Lambda^s = \Cal L_s $ becomes
$$ 
\Cal L_s  = \bar\ell_s\ell_s
+ (\be_1+A)\bar\ell_s +(\be_2+A)\ell_s + (\be_3+A)
$$
where $\bar\ell_s,\ell_s$ are first-order differential
operators differing from $\bar\ell,\ell$ respectively
by terms of order zero, and taking the forms
$\bar\ell_s = \p_x + i\al(t\p_t-s)$,
$\ell_s = \p_x -i\bar\al(t\p_t-s)$
for $|x|\le r/2$,
where $\al$ depends only on $t$
and the $\be_i$ only on $x$.

To see that $\bar\ell_s$ does take the form claimed
for $|x|\le r/2$, express
$\bar\ell = \p_x +i\al t\p_t$ modulo terms
$\gamma(x,t)\p_x$ and $\gamma(x,t)\p_t$ where
$\gamma\equiv 0$ for $|x|\le r/2$. 
Then $\Lambda^{-s}[\gamma(x,t)\p_x, \Lambda^s]$ 
and $\Lambda^{-s}[\gamma(x,t)\p_t, \Lambda^s]$
are operators of the type $A$, since they have
nonpositive orders and their symbols of order
zero vanish identically for $|x|<r/2$.
\hfill\hfill\qed

\subheading{4. A Two-Dimensional Problem And Preliminary Inequalities}

The remainder of the paper consists of a self-contained
analysis of a special class of pseudodifferential equations
in a real two-dimensional region. We begin by describing
the equations in question and fixing notation, which in some
respects differs from that of preceding sections. 

Fix an interval $I=[-r,r] \subset\rr$.
Denote by $(x,t)\in\rr^2$ coordinates in a neighborhood
$U$ of $\Cal D = I\times\{0\}$. 
The interval $\Cal D$ corresponds to the degenerate annulus 
embedded in the boundary of the worm domain, and will be the
focus of attention.
The convention concerning the symbols $A,B,E$ introduced 
at the outset of \S3 remains in force.

Consider a one parameter family of pseudodifferential operators 
of the form
$$
\Cal L_s = \bar L L
+(\be_1(x)+A)\bar L + (\be_2(x)+A) L 
+ (\be_3(x)+A)
\tag4.1
$$
where the $\be_j$ are $C^\nf$ functions.
Suppose that 
$\bar L,\ L$ are first-order differential operators
depending on the real parameter $s$, and that $-L$ is
the formal adjoint of $\bar L$, modulo an operator of order zero.
Suppose that where $|x|\le r$, 
they take the special forms\footnote{No assumption is now made
on the vanishing or nonvanishing 
of the coefficient of $\p_t$ in $\bar L$ where $|x|>r$, but
the strict pseudoconvexity of $\worm$ outside the exceptional annulus 
was used to reduce Proposition 1 to Proposition 2 below.}
$$
\align
\bar L &=
\p_x + ia(x)(t\p_t+s) + O(t^2)\p_t   \\
L &=
\p_x - i\bar a(x)(t\p_t+s) + O(t^2)\p_t.
\endalign
$$
Here $O(t^2)$ denotes multiplication by a smooth function 
divisible by $t^2$ on the region $U$.
$a$ and the coefficients $\be_j$ are assumed independent of $s$, 
but $A$ and the terms $O(t^2)\p_t$ are permitted to depend on $s$.

Assume that 
$$
\Re a(x)\ne 0 \qquad \text{for all } x\in I, \tag4.2
$$
and that there exist smooth real-valued coefficients $\mu,\nu$ 
such that $[\bar L,L] = i\mu(x,t)\p_t + i\nu(x,t)\Re\bar L$,
satisfying
$$
\mu\ge 0 \text{ at every point of } U.  \tag4.3
$$

Because $\bar L^* = - L$ modulo
a term of order zero,
$L$ has the same real part as $\bar L$.
A change of variables of the form
$(x,t)\mapsto (x,h(x,t))$, with $h(x,0)\equiv 0$ where
$|x|\le r$, therefore
reduces matters to the case where the real parts of both $\bar L$
and $L$ are everywhere parallel to $\p_x$,
and $\bar L = \p_x +i\tilde a(x,t)(t\p_t+s) + O(t^2)\p_t$ 
on $I\times\rr$,
with $\tilde a$ real-valued and nonvanishing.
Rewrite $\tilde a(x,t) =  a(x) + O(t)$, and incorporate
the contribution of $O(t)$ into the various terms
$O(t^2)\p_t$ and $A$. (4.3) is invariant under diffeomorphism
and hence $\mu$ cannot change sign, so the coefficient
of $t$ in the Taylor expansion of $\mu(x,t)$ about $t=0$
must vanish identically, for $|x|\le r$.
This forces $\p_xa(x)\equiv 0$ there. Thus
$$
\bar L = \p_x +i a(t\p_t+s) + O(t^2)\p_t  \tag4.4
$$
for $|x|\le r$, where $a$ is a nonzero real constant. 
Moreover, $\p_x$ may be expressed in $U$ as a nonvanishing
scalar multiple of $\bar L + L$, modulo an operator of order zero.
>From now on we work in these new coordinates.

Define $\Gamma=\{(x,t,\xi,\tau): (x,t)\in\Cal D \text{ and }
\xi=0\}$.
Decompose $\Gamma=\Gamma^+\cup\Gamma^-$ where $\Gamma^+
=\{\tau>0\}\cap\Gamma$.
Then by (4.3), in some conic neighborhood of $\Gamma^+$
the principal symbol of $[\bar L,\,L]$ 
equals a nonpositive symbol, modulo terms in the span
of the symbols of $\bar L,L$ and a term of order zero.

The symbol $\|\cdot\|$, with no subscript, denotes the norm
in $L^2(U)$, while $\|\cdot\|_t$ denotes any fixed
norm for the Sobolev space $H^t$ of functions having $t$
derivatives in $L^2$ and supported in $U$.
The goal of the remainder of the paper
is the following {\it a priori} estimate.

\proclaim{Proposition 2}
Let $\{\Cal L_s\}$ be a family of operators of the form (4.1)
satisfying all of the hypotheses introduced above.
Then there exist a discrete exceptional set $S\subset[0,\nf)$ 
such that for any $s\notin S$ 
and any pseudodifferential operator $Q$ 
of order $0$ whose principal symbol is nonzero
in some conic neighborhood of $\Gamma^-$, 
there exist $C<\nf$ and a neighborhood $W$ of $\Cal D$
such that for every $C^\nf$ function $u$ supported in $W$,
$$
\|u\| + \|\bar L u\| + \|Lu\|
\le C\cdot \left( \|\Cal L_s u\| + \|u\|_{-1} 
+ \|Qu\|_1 \right).
$$
\endproclaim

The key conclusions are that there is no loss of derivatives
in estimating $u$ in terms of $\Cal L_s u$, and that this
holds for a sequence of values of $s$ tending to $+\nf$.
The assumption that $u\in C^\nf$ is essential.
All hypotheses of \S4 are satisfied by the family of operators
$\Cal L_s$ derived in \S 2 and \S 3.
Proposition 2 thus implies the validity of (3.2), and hence
of Proposition 1, which in turn implies our Theorem.

Our first preliminary estimate is a standard one valid for all $s\in\rr$.

\proclaim{Lemma 1} For each exponent $s$ and each $Q$
there exists $C<\nf$ such that
$$
\|\p_x u\|
\le C\|u\| + C\|\Cal L_s u\| + C \|Qu\|_1
$$
for every $u\in C_0^\nf(U)$.
\endproclaim

\demo{Proof} 
For $(x,t)\in\Cal D$,
$\sigma_2(\Cal L_s)(x,t,\xi,\tau)=0$ if and only if
$(x,t,\xi,\tau)\in\Gamma$. Therefore
the characteristic variety of $\Cal L_s$ in
$T^*W$ is contained in an arbitrarily small conic neighborhood
of $\Gamma$ as $\de\to 0$. 
Consequently there exists an operator $\tilde Q$ of order zero 
such that firstly, $T^*W$ is contained in the union of the
two regions where $\tilde Q$ is elliptic and $\tau>0$,
and secondly, the symbol of $\tilde Q$ is supported in the union 
of the two regions where $\Cal L_s$ is elliptic, 
and where $Q$ is elliptic.

Since $\Cal L_s$ is elliptic outside a small conic neighborhood
of $\Gamma$, the $H^2$ norm of $u$ is 
majorized away from $\Gamma$ by $\|\Cal L_s u\| + \|u\|_{-1}$,
while in a conic neighborhood of $\Gamma^-$ the $H^1$ norm of $u$ is  
majorized by $\|Qu\|_1+ \|u\|_{-1}$. 

Write $\lan f,\,g\ra = \int_{U} f \bar g\,dx\,dt$.
By G\aa rding's inequality and the fact that $i\mu\cdot i\tau
\le 0$ in the support of the symbol of $\tilde Q$,
$$
- \Re\big( \lan \bar L Lu,\,u \ra \big)
\ge c\|Lu\|^2 + c\|\bar Lu\|^2 -C\|\bar Lu\|\cdot\|u\|
-C\|Lu\|\cdot\|u\|
-C\|u\|^2 - C\|\tilde Qu\|_1^2.
$$
The second condition imposed on $\tilde Q$ ensures that
$$\|\tilde Qu\|_1
\le C\|\Cal L_s u\| + C\|Qu\|_1 + C \| u\|_{-1}.$$
Estimating $\lan (\Cal L_s - \bar L L)u,\,u\ra$ by
Cauchy-Schwarz thus leads to
$$\|\bar L u\| + \|L u\| \le C\|\Cal L_s u\| + C\|u\|
+ C\|Qu\|_1.$$
But $\p_x$ may be expressed as a linear combination of $\bar L$
and of $L$ modulo an operator of order $0$.
\hfill\hfill\qed\enddemo

\proclaim{Lemma 2}
There exists $C<\nf$ such that for every $f\in C^1(\rr)$
and every $\e>0$,
$$
\|f\|_{\lt[\e,2\e]}
\le C\|f\|_{\lt[-2\e,-\e]} + C \e \|\p_x f\|_{\lt(\rr)}.
$$
Likewise
$$|f(0)-f(-\e)| \le C\e^{1/2}\|\p_x f\|_{\lt}.$$
\endproclaim

The conclusions are invariant under translation, and the lemma
will be invoked in that more general form. 

\demo{Proof} 
For each $x\in (\e,2\e)$,
$|f(x) - f(x-3\e)|\le\int_{-2\e}^{2\e} |\p_x f(y)|\,dy$
and both conclusions follow from the triangle and Cauchy-Schwarz
inequalities.
\hfill\hfill\qed \enddemo

To simplify notation define
$$
\okay = \|\Cal L_s u\| + \|u\|_{-1} + \|Qu\|_1.
$$
Let $\de>0$ be a small constant to be chosen in \S6,
and assume $u$ to be supported in
$$  W\subset\{(x,t): |t|<\de,\ |x|<r+\de\}. $$

Applying Lemma 2 to the function $x\mapsto u(x,t)$ for each $t$
and applying Lemma 1 gives the following estimate,
under the hypotheses of Lemma 1.

\proclaim{Lemma 3} 
$$
\|\p_x u\|  + \|u\| \le C\|u\|_{\lt(I\times(-\de,\de))} + C\okay.
$$
\endproclaim

\subheading{5. Limiting Operators and Mellin Transform}

Let $a$ be a nonvanishing $C^\infty$, real-valued function.
For $\z\in\CC$ define the ordinary differential operator
$$
H_\z = 
(\p_x + i\z a(x)) (\p_x-i\z a(x)) +\be_1(x) (\p_x + i\z a(x))
+\be_2(x) (\p_x-i\z a(x)) +\be_3(x),
$$
acting on functions of $x\in I$. 
Only the case of constant $a$ will be needed in this
paper, but the general case arises in another problem and
hence merits discussion.

\proclaim{Definition 1}
$\frak S$ is defined to be the set of all $\z\in\CC$ such that
there exists a solution $g$ of $H_\z g\equiv 0$ on $I$,
satisfying $g(-r)=g(r)=0$.
\endproclaim

For any complex number $w$ we write
$\jj w. = (1+|w|^2)^\half$.

\proclaim{Lemma 4}
$\frak S$ is a discrete subset of $\CC$, and for any compact
subset $K$ of $[0,\nf)$, the set of all $\z\in\frak S$ having
real part in $K$ is finite. 
For each $s$ such that
$\frak S\cap (s+i\rr) = \emptyset$ there exists
$C<\nf$ such that 
for all $\z \in s+i\rr$,
for all $f,\varphi,\psi\in C^\nf(I)$ 
satisfying $H_\z f = \varphi + \p_x\psi$,
one has
$$\aligned
&\|f\|_{\lt(I)} 
+ \jz\rp \|\p_x f\|_{\lt(I)} \\
&\qquad\qquad 
\le C\jz^{-1/2} \big(|f(-r)| +  |f(r)|\big)
+ C\jz^{-2}\|\varphi\|_{\lt(I)}
+ C\jz^{-1}\|\psi\|_{\lt(I)}. 
\endaligned
\tag5.1
$$
\endproclaim

\demo{Proof}
Throughout this proof, all norms without subscripts
denote $\lt$ norms.
The selfadjoint part of $H_\z$, applied to $f$, equals 
$(\p_x-\g a(x))(\p_x+\g a(x))f$, modulo $O(\jj \g. \|f\| + \|\p_x f\|)$,
in the $\lt(I)$ norm.
Therefore for $s$ in any fixed compact subset of $\rr$ and
any $\z= s+i \g \in s+i\rr$, 
for any $f$ vanishing at both endpoints of $I$,
$$
-\Re\lan H_\z f,\,f\ra
\ge \|\p_x f\|^2 + \g^2 \int_I |f|^2|a|^2
-O(\jj \g.\|f\|^2 + \|f\|\cdot\|\p_x f\|).
$$
The coefficient $a$ vanishes nowhere,
while
$$|\lan H_\z f,\,f\ra| = |\lan \varphi+\p_x\psi,\,f\ra|
\le \|f\|\cdot\|\varphi\| +
\|\p_x f\|\cdot\|\psi\|.$$
Combining the last two inequalities and invoking the
Cauchy-Schwarz inequality and small constant -- large constant
trick, one obtains
$$
\g^2\|f\|+ |\g|\cdot \|\p_x f\| \le C\|\varphi\| + C|\g|\cdot\|\psi\| 
\tag5.2
$$
for all sufficiently large $|\g|$, under the additional hypothesis
that $f$ vanishes at both endpoints of $I$.

There exists a unique solution $\phi_\z$ of $H_\z \phi_\z=0$
on $I$, satisfying $\phi_\z(-r)=0$, $\p_x\phi_\z(-r)=1$. Then
$\phi_\z(r)$ is an entire holomorphic function of $\z$, 
and $\z\in \frak S \Leftrightarrow \phi_\z(r)=0$.
We have seen that $\z\notin\frak S$ provided that the imaginary
part of $\z$ is sufficiently large, when the real part stays
in a bounded set. Thus $\phi_\z(r)$ is nonconstant, so has
discrete zeros.\footnote{An alternative
method of proof would be to combine (5.2) with general results
from the perturbation theory of linear operators \cite{Ka}, utilizing
again the holomorphic dependence of $H_\z$ on $\z$.}

To prove (5.1) let $\z=s+i\g$ and $f\in C^2(I)$
be given, and decompose $f = g + h$ where 
$H_\z g\equiv 0$ and $h$ vanishes at the endpoints of 
$I$. The hypothesis $\frak S\cap (s+i\rr)=\emptyset$ means
that the Dirichlet nullspace of $H_\z$ is $\{0\}$,
so by elementary reasoning we conclude that
for each $\g$ there exists $C<\nf$ such that 
$\|h\|+ \|\p_x h\|\le C\|\varphi\| + \|\psi\|$,
since $H_\z h = H_\z f = \varphi + \p_x\psi$. 
Moreover, since $H_\z$ depends continuously
on $\z$, $C$ may be taken to be independent of $\z$ 
in any compact subset of $\CC\backslash \frak S$.
When $|\g|$ is sufficiently large, on the other hand, 
(5.2) implies
$\|h\| + \jz\rp\|\p_x h\|\le C\jz^{-2} \|\varphi\| + C\jz^{-1}\|\psi\|$. 
Thus the component $h$ of $f$ satisfies (5.1). 

We have $H_\z g=0$, so that clearly $\|g\|$ and $\|\p_x g\|$
are majorized by $C|f(r)| + C|f(-r)|$,
uniformly for $\z$ in any compact set disjoint from $\frak S$.
Assuming henceforth that $|\g|$ is large, 
the equation gives the inequality
$$\|\p_x^2 g\|\le C\g^2\|g\| + C|\g| \cdot \|\p_x g\|.$$
Integrating by parts as in the proof of (5.2) yields
$$
\|\p_x g\|^2 + \g^2\|g\|^2 \le C|g(-r)\p_x g(-r)| + 
C|g(r)\p_x g(r)|. \tag5.3
$$

To control the right hand side we use the bound
$$
|\p_x g(-r)| + |\p_x g(r)|
\le C|\g|^{1/2}\|\p_x g\| + C |\g|^{-1/2}\|\p_x^2 g\|.
$$
Indeed, setting $v=\p_x g$, for any $r'\in[r-|\g|^{-1},r]$
$$
|v(r)-v(r')|
\le C\int_{r'}^r |\p_x v| 
\le C|\g|^{-1/2} \|\p_x v\|_{L^2}.
$$
Then
$$
|v(r)|
  \le |\g|\int_{r-|\g|^{-1}}^r |v(r)-v(r')| \,dr'
  + |\g|\int_{r-|\g|^{-1}}^r |v(r')|\,dr' 
$$
and the desired bound follows by Cauchy-Schwarz.

Putting this into (5.3), introducing a parameter $\la\in\rr^+$
and applying Cauchy-Schwarz yields
$$\align
\g^2\|g\|^2 &+ \|\p_x g\|^2 \\
&\le C\la|g(-r)|^2 + C\la|g(r)|^2
+ C\la\rp|\g| \cdot \|\p_x g\|^2 + C\la\rp|\g|\rp\|\p_x^2 g\|^2  \\
& \le  
C\la |g(-r)|^2 + C\la|g(r)|^2 
+ C\la\rp |\g| \cdot \|\p_x g\|^2
+ C\la\rp|\g|\rp \big(\g^4\|g\|^2 + \g^2\|\p_x g\|^2\big) \\
& \le
C\la |g(-r)|^2 + C\la|g(r)|^2 
+ C\la\rp |\g| \cdot \|\p_x g\|^2 + C\la\rp|\g|^3 \|g\|^2.
\endalign
$$
Choose $\la$ to be a large constant times $|\g|$. Then
the last two terms on the right-hand side may be absorbed into
the left, leaving
$$
\g^2\|g\|^2 + \|\p_x g\|^2 \le C|\g| \cdot |g(-r)|^2
+ C|\g| \cdot |g(r)|^2.  
$$
Since $g=f$ at the endpoints of $I$, this is the desired inequality
for $g$. Adding it to that for $h$ concludes the proof.
\hfill\hfill\qed\enddemo

\proclaim{Definition 2}
$$
S=\{s\in[0,\nf): \text{ there exists } \g\in\rr \text{ such that }
s-\tfrac12 + i\g\in \frak S\}.
$$
\endproclaim
\noindent Lemma 4 guarantees that $S$ is discrete.

Specialize now to the case where $a(x)\equiv a$, the
real constant in (4.4). Define
$$\align
\bl_s  =
& (\p_x + ia(t\p_t+s))\circ(\p_x-ia(t\p_t+s)) \\
& \qquad +\be_1(x)(\p_x + ia(t\p_t+s))
+\be_2(x)(\p_x-ia(t\p_t+s))
+\be_3(x).
\endalign $$
Expanding the last term in the expression 
$\bl_s = \Cal L_s + (\bl_s - \Cal L_s)$ gives
$$ \bl_s u = \Phi + \p_x\Psi  $$
where
$$\aligned
\Phi &= \Cal L_s u + (t\p_t)^2 A u + t\p_t A u + A u  \\
\Psi & =  t\p_t A u +  A u
\endaligned
\tag5.4
$$
To reach (5.4) we may for instance express
$t\p_t\circ O(t^2)\p_t$ as $(t\p_t)^2A + t\p_tA + A$,
since multiplication by $t$ is an operator of the type $A$. 
Likewise $[A, t\p_t] = t[A,\p_t] + [A,t]\p_t$ is an operator 
of type $A$, because $\sigma_{-1}([A,t]) = c\p_{\tau}\sigma_0(A)$
vanishes identically for $(x,t)\in\Cal D$ since 
$\sigma_0(A)$ itself vanishes there.

The partial Mellin transform of $f$ with
respect to the $t$ variable is defined to be
$$\hat f(x,\g) = \int_0^\nf f(x,t)t^{-i\g}\,\dt,$$
provided that the integral converges.
If $f(x,\cdot)\in C^\nf[0,\nf)$ has bounded support for each $x$, 
then the integral defining $\hat f(x,\g)$ converges absolutely
whenever $\g$ has strictly positive imaginary part, and
$\hat f(x,\g)$ extends to a meromorphic function of $\g\in\CC$,
whose only possible poles are at $\g=0,-i,-2i,\dots$. 
Clearly
$$(t\p_t f)\sphat\,(x,\g) = i\g \hat f(x,\g)$$
for all such $f$.  Consequently
$$(\bl_s u)\sphat\,(x,\g) = H_{s+i\g}\hat u(x,\g)
\qquad\text{ for all } \g\in\CC\backslash\{0,-i,-2i,\dots\}.$$

The Mellin inversion and Plancherel formulas read
$$
f(x,t) = c\int_\rr \hat f(x,\g)\,t^{i\g}\,d\g, \qquad
\int_0^\nf |f(x,t)|^2\,\dt 
= c'\int_\rr |\hat f(x,\g)|^2\,d\g.
$$
It follows directly from the definitions that
$(t^\half f)\sphat\,(x,\g) = \hat f(x,\g+\tfrac{i}2)$
for all $\g\in\rr$.
Thus the Plancherel identity may be rewritten as
$$
\int_0^\nf |f(x,t)|^2\,dt
= \int_0^\nf |t^\half f(x,t)|^2\,\dt 
= c'\int_\rr |\hat f(x,\g+\tfrac{i}2)|^2\,d\g.
$$

\subheading{6. Proof of the Main Estimate}

We may now estimate $u$ in terms of $\Cal L_s u$.
To begin,
$$
\iint_{I\times[0,\de)} |u(x,t)|^2\,dx\,dt
= c'\iint_{I\times\rr} |\hat u(x,\g+\tfrac{i}2)|^2\,d\g\,dx.
$$
Assume that $s\notin S$, and write $\z=s-\tfrac12+i\g$.
With $\Phi,\Psi$ defined as in (5.4),
$H_\z\hat u(x,\g+\tfrac{i}2) = \hat\Phi(x,\g+\tfrac{i}2)
+ \p_x\hat\Psi(x,\g+\tfrac{i}2)$.
Applying Lemma 4 on $I$ yields for each $\g\in\rr$
$$\aligned
\int_{I}\hat u(x,\g+\tfrac{i}2)\,dx
& \le C\int_{I} 
|\hat\Phi(x,\g+\tfrac{i}2)|^2\jj \g.^{-4} \,dx  
 + C \int_{I} 
|\hat\Psi(x,\g+\tfrac{i}2)|^2\jj \g.^{-2} \,dx \\
&\qquad + C |\hat u(-r, \g+\tfrac{i}2)|^2
+ C |\hat u(r, \g+\tfrac{i}2)|^2.
\endaligned
\tag6.1
$$

\proclaim{Lemma 5}
Assume $W\subset\{(x,t): |t|<\de,\ |x| < r+\de \}$.
Then there exists $C<\nf$ such that 
$$\iint_{I\times\rr}
|\hat\Phi(x,\g+\tfrac{i}2)|^2\jj \g.^{-4}
\,d\g\,dx
\le C \de^2\|u\|^2 + C\okay^2
$$
and
$$\iint_{I\times\rr}
 |\hat\Psi(x,\g+\tfrac{i}2)|^2\jj \g.^{-2}
\,d\g\,dx
\le C \de^2 \|\p_x u\|^2 + C\okay^2.
$$
\endproclaim

Granting the lemma, we conclude from (6.1) that
$$
\|u\|_{\lt(I\times[0,\de))}^2
\le C \de^2 \|u\|^2 + C \de^2 \|\p_x u\|^2
+ C\okay^2
+  C\int_{\rr} |u( r,t)|^2\,dt
+  C\int_{\rr} |u(-r,t)|^2\,dt.
$$
But by Lemma 2 and the assumption that $u(x,t)\equiv 0$ for
$|x|>r+\de$, this last term is dominated by $C\de \|\p_x u\|^2$.
Thus
$$
\|u\|_{\lt(I\times[0,\de))}^2
\le C \de^2 \|u\|^2 + C \de\|\p_x u\|^2
+ C\okay^2.
$$
All the same reasoning applies on the region $I\times (-\de,0]$,
after the change of variables $t\mapsto -t$. Thus
$$
\|u\|_{\lt(I\times(-\de,\de))}
\le C\de^\half (\|u\| + \|\p_x u\|)
+ C\okay.
$$
Combining this with Lemma 3 gives
$$\|u\|+ \|\p_x u\|
\le C \de^\half (\|u\| + \|\p_x u\|)
+ C\okay,
$$
so choosing $\de$ to be sufficiently small 
gives $\|u\|\le C \okay$, concluding the proof.
\hfill\hfill\qed\medskip

\demo{Proof of Lemma 5}
The principal term in the double integral of the lemma for
$\Phi$ is of course the contribution of $\Cal L_s u$:
$$\align
&\iint_{I\times\rr} \left|
(\Cal L_s u)\sphat\,(x,\g+\tfrac{i}2) \right|^2
\jj \g.^{-4}\,d\g\,dx \\
&\qquad\qquad \le 
\iint_{I\times\rr} \left|
(\Cal L_s u)\sphat\,(x,\g+\tfrac{i}2) \right|^2
\,d\g\,dx \quad
\le \quad C\|\Cal L_s u\|^2,
\endalign
$$
as desired.

Any operator $A$ of order $\le 0$ satisfying 
$\sigma_0(A)(x,t,\xi,\tau)\equiv 0$ for $(x,t)\in\Cal D$ satisfies
$$
\|Au\|^2 \le C\de^2 \|u\|^2 + C\|u\|_{-1}^2
$$
for all $u$ supported in $W$, as $\de\to 0$.
A typical term of $\Phi$ resulting from $\bl_s u-\Cal L_s u$
is $(t\p_t)^2 Au$. Its contribution to the first double
integral in Lemma 5 is
$$\align
\iint_{I\times\rr} \left|
((t\p_t)^2 Au)\sphat\,(x,\g+\tfrac{i}2) \right|^2
\jj \g.^{-4}\,d\g\,dx  
&\le
C \iint_{I\times\rr} \left|
(Au)\sphat\,(x,\g+\tfrac{i}2) \right|^2 \,d\g\,dx \\
& 
= C\iint_{I\times [0,\de)} |Au(x,t)|^2 \,dx\,dt \\
&
\le C \de^2 \| u\|^2 + C\|u\|_{-1}^2. 
\endalign $$
A typical constituent of $\Psi$ is the term
$t\p_t Au$. Its contribution is dominated by
$$
C \iint_{I\times\rr} \left|
(t\p_t Au)\sphat\,(x,\g+\tfrac{i}2) \right|^2
\jj \g.^{-2}\, d\g\,dx 
\le
C \iint_{I\times\rr} \left|
(Au)\sphat\,(x,\g+\tfrac{i}2) \right|^2 \,d\g\,dx   
$$
and the remainder of the calculation is as above.
\hfill\hfill\qed\enddemo

\bigskip
\noindent{\smc Comments.}
The author can advance no reason why the method of reduction to
the boundary should be essential to this analysis. Working
directly on $\overline{\worm}$ might well result in a shorter proof.
On the other hand, the analysis in \S\S4-6 applies with minor
modification to a broader class of equations unconnected with
the $\bar\p$--Neumann problem.

Some refinements of this analysis and related observations
will appear in \cite{Ch2}.

\Refs

\k B. \by D.~Barrett
\paper Behavior of the Bergman projection on the 
Diederich-Forn\ae ss worm
\jour Acta Math.\vol168\yr1992\pages1-10 \er


\k BL. \by S.~Bell and E.~Ligocka
\paper A simplification and extension of Fefferman's theorems
on biholomorphic mappings
\jour Invent. Math.\vol57\yr1980\pages283-289 \er

\k BS1. \by H.~Boas and E.~Straube
\paper Equivalence of regularity for the Bergman projection and
the $\bar\p$--Neumann operator 
\jour Manuscripta Math.\vol67\yr1990\pages25-33 \er

\k BS2. \bysame
\paper Sobolev estimates for the $\bar\p$--Neumann operator
on domains in $\bold C^n$ admitting a defining function
that is plurisubharmonic on the boundary
\jour Math. Zeitschrift\vol206\yr1991\pages81-88 \er

\k Ca1. \by D.~Catlin
\paper Subelliptic estimates for the $\bar\p$--Neumann problem
on pseudoconvex domains
\jour Annals of Math.\vol126\yr1987\pages131-191 \er

\k Ca2. \bysame
\paper Global regularity of the $\bar\p$--Neumann problem
\jour Proc. Symp. Pure Math.\vol41\yr1984\pages39-49 \er

\k CNS.  \by D.-C.~Chang, A.~Nagel and E.~M.~Stein
\paper Estimates for the $\bar\partial$--Neumann problem in
pseudoconvex domains of finite type in $\CC^2$
\jour Acta Math. \vol169\yr1992\pages153-228 \er

\k Ch1. \by M.~Christ
\paper The Szeg\"o projection need not preserve global
analyticity \jour Annals of Math. \toappear \er

\k Ch2. \bysame
\paper Global $C^\nf$ irregularity for mildly degenerate
elliptic operators \paperinfo preprint
\er

\k DF.  \by K.~Diederich and J.~E.~Forn\ae ss
\paper Pseudoconvex domains: an example with nontrivial
Nebenh\"ulle
\jour Math. Ann.\vol225\yr1977\pages275-292 \er

\k FS. \by J.~E.~Forn\ae ss and B.~Stens\o nes
\book Lectures On Counterexamples In Several Complex
Variables
\publ Princeton University Press\publaddr Princeton, NJ
\bookinfo Mathematical Notes 33 \yr1987
\er

\k Ka. \by T.~Kato 
\book Perturbation Theory For Linear Operators 
\publ Springer-Verlag\publaddr New York \yr1966 \er

\k K1. \by J.~J.~Kohn
\paper Global regularity for $\bar\p$ on weakly
pseudoconvex manifolds
\jour Trans. AMS\vol181\yr1973\pages273-292 \er

\k K2.  \bysame
\paper Estimates for $\bar\p_b$ on pseudoconvex CR manifolds
\jour Proc. Symp. Pure Math.\vol43\yr1985\pages207-217 \er


\endRefs
\enddocument\end